# Integral points on hyperbolas over $Z$ :
# A special case

*Konstantine Zelator*
*Department of Mathematics and Computer Science*
*Rhode Island College*
*600 Mount Pleasant Avenue*
*Providence, R.I. 02908-1991, U.S.A.*
*e-mail address: 1) Kzelator@ric.edu*
*2) Konstantine_zelator@yahoo.com*

*Effective August 1, 2009, and for the academic year 2009-10;*
*Konstantine Zelator*
*Department of Mathematics*
*301 Thackeray Hall*
*139 University Place*
*University of Pittsburgh*
*Pittsburgh, PA 15260*
*U.S.A*
*e-mail address : kzet159@pitt.edu*

1. **Introduction**

   The subject matter of this work is integral point on conics described by the standard general form equation

   $$\alpha x^2 + \beta xy + \gamma y^2 + \delta x + \varepsilon y + J = 0 \tag{1}$$

   and where the coefficients $\alpha, \beta, \gamma, \delta, \varepsilon, J$ are integers satisfying the conditions

   $$\begin{cases} \beta^2 - 4\alpha\gamma = k^2, \\ \text{with } \alpha \neq 0, \gamma \neq 0, \text{and } k \text{ a positive integer} \end{cases} \tag{2}$$

   It is a well-known fact if a conic described by (1), with the six coefficients being real numbers and $\beta^2 - 4\alpha\gamma > 0$; then such a conic must be either a hyperbola or a pair of two intersecting straight lines (the degenerate case). There is an extensive body of literature on conics, spanning a few hundred years; but with most books on the subject having been published in the last 150 years. For example, the reader may refer to **[1]**. An *integral point* on a conic is simply an ordered pair $(x, y)$ satisfying (1) and with both $x$ and $y$ being integers. Note that because of condition (2), equation (1) when considered as a diophantine equation in the variables $x$ and $y$, cannot be a *Pell* equation; since in the case of Pell equation, $\alpha = 1, \gamma = -d$, $\beta = \delta = \varepsilon = 0$, and $J = 1$ or -1; where $d$ is a positive integer which is not a perfect or integer square; and so (2) could not be satisfied. A hyperbola which is described by a Pell-equation of the form $x^2 - dy^2 = 1$ ($d$ a nonperfect square), has in fact infinitely integral points. A proof of this fact, as well as the method of finding all the solutions, can be found in a number of number theory books. As examples, see [2] or [3]. On the other hand, a Pell equation of the form $x^2 - dy^2 = -1$ will either have infinitely many integer solutions or no solutions at all; depending on the period in the simple continued fraction expansion of the irrational number $\sqrt{d}$ (see references [2] and [3]). Again, as mentioned before, none of the hyperbolas we study in this paper corresponds to a Pell equation. As we will see, the integer $I = k^2 \cdot (\delta^2 - 4\alpha J) - (2\alpha\varepsilon - \beta\delta)^2$, plays a key role. When $I \neq 0$, each of hyperbolas described by (1), and under the conditions in (2), has finitely many integral points (including the case of zero or no integral points). These finitely many integer pairs $(x, y)$ can be found by a method of technique outlined in Sections 2 and 3. This method uses only straightforward algebra and a couple of basic facts on quadratic trinomials. In Section 4 we offer a numerical example. In Section 5, we offer a couple of observations and remarks on ther integer $I$ above. In Section 7 we take a look at the special case $\alpha = 1 = k$. And in Section 8, we consider the case $I = 0$. When $I = 0$, a conic described by (1) and (2), now becomes a pair or union of two intersecting straight lines. A straight line is described by an equation of the form $ax + by = c$; and in our case with $a, b, c$ being integers. Such a line can either have no or infinitely many integral points. Solving such a linear diophantine equation is standard, well known material which can be found in every introductory number theory book. So in the case $I = 0$, the curve in question is the union of two intersecting straight lines.

$$\left. \begin{array}{l} l_1 : a_1 x + b_1 y = c_1 \\ l_2 : a_2 x + b_2 y = c_2. \end{array} \right\} ; \text{with } a_1, b_1, a_2, b_2, c_1, c_2 \text{ in } Z$$

These are exactly four possibilities: $l_1$ contains infinitely many integral points and $l_2$ none; or vice-versa; or each of them contains infinitely many integral points; or neither of them does.

2. **A solving technique**
   Equation (1) is equivalent to

$$\alpha x^2 + \beta xy + \gamma y^2 + \delta x + \varepsilon y + J + \lambda = \lambda \qquad (3)$$

where $\lambda$ is a rational number to be determined.
Consider the left-hand side of (3) as a quadratic trinomial in $x$.
We write it in standard form:

$$\alpha x^2 + (\beta y + \delta) x + (\gamma y^2 + \varepsilon y + J + \lambda) = \lambda \qquad (4)$$

The discriminant $D(y)$ of this trinomial in $x$, depends on $y$.
We have $D(y) = (\beta y + \delta)^2 - 4\alpha(\gamma y^2 + \varepsilon y + J + \lambda)$, and by (2) we obtain

$$D(y) = k^2 y^2 - 2(2\alpha\varepsilon - \beta\delta) y + (\delta^2 - 4\alpha J - 4\alpha\lambda) \qquad (5)$$

The idea here is pretty simple: to choose a rational number $\lambda$ such that $D(y)$ becomes the square of a linear polynomial in $y$ and of the form $ky + r$, where $r$ is a rational number; so that we may further factor the trinomial on the left-hand side of (4) as a product of two linear polynomials in $x$ and $y$ with integer coefficients.
We choose $\lambda$ such that

$$\delta^2 - 4\alpha J - 4\alpha\lambda = \left( \frac{2\alpha\varepsilon - \beta\delta}{k} \right)^2 \qquad (6)$$

Then from (5) we obtain,

$$D(y) = \left[ ky - \left( \frac{2\alpha\varepsilon - \beta\delta}{k} \right) \right]^2 \qquad (7)$$

Solving for $\lambda$ in (6) produces,

$$\lambda = \frac{k^2(\delta^2 - 4\alpha J) - (2\alpha\varepsilon - \beta\delta)^2}{4\alpha k^2} \quad (8)$$

The left-hand side of equation (4) can be factored, according to the fundamental Theorem of Algebra, as

$$\alpha(x - r_1(y))(x - r_2(y)) = \lambda \quad (9)$$

where $r_1(y) = \dfrac{-(\beta y + \delta) + \sqrt{D(y)}}{2\alpha}, r_2(y) = \dfrac{-(\beta y + \delta) - \sqrt{D(y)}}{2\alpha};$

are the two roots (which of course depend on $y$) of the quadratic trinomial in $x$. By (9),(8),(7), and (4); and after some algebra (which includes multiplying both sides of (9) by $4\alpha k^2$) we arrive at the desired factorization:

$$[2\alpha k + k(\beta - k)y + \delta k + 2\alpha\varepsilon - \beta\delta] \cdot [2\alpha kx + k(\beta + k)y + \delta k - (2\alpha\varepsilon - \beta\delta)]$$
$$= k^2(\delta^2 - 4\alpha J) - (2\alpha\varepsilon - \beta\delta)^2 \quad (10)$$

Let $I = k^2(\delta^2 - 4\alpha J) - (2\alpha\varepsilon - \beta\delta)^2 \quad (10a)$

Note that if one were to start with equation (10); and multiply out the two factors of the left-hand side; use (2) and collect the aggregate $x^2, y^2, xy, x, y,$ and constant terms; each of these terms will contain a common factor $4\alpha k^2$; after the cancellation of which, one returns to equation (1).

3. **The case** $I \neq 0$

When $I = k^2(\delta^2 - 4\alpha J) - (2\alpha\varepsilon - \beta\delta)^2 \neq 0$, and $(x, y)$ is an integer solution of (1), and hence of (10) as well (and conversely), then the two factors on the left-hand side of equation (10); will be nonzero divisors of $I$, whose product is $I$.

Let $1 = d_1 < d_2 < \ldots < d_N = |I|$ be the $N$ positive integers which are the positive divisors of the natural number $|I|$. In order to find all the integer pairs $(x, y)$ satisfying (10), we must solve $2N$ linear systems in the unknowns or variables $x$ and $y$. We group these $2N$ systems into $N$ groups, each group containing two systems:

The ith group $\begin{cases} 2\alpha kx + k(\beta - k)y + \delta k + 2\alpha\varepsilon - \beta = e \cdot d_i \\ 2\alpha kx + k(\beta + k)y + \delta k - (2\alpha\varepsilon - \beta) = e \cdot \dfrac{I}{d_i} \end{cases},$

Where $e = 1$ or $-1$ (which explains why each group contains two systems); or equivalently, in standard form,

$$\begin{cases} 2\alpha kx + k(\beta - k)y = e \cdot d_i - \delta k - (2\alpha\varepsilon - \beta\delta) \\ 2\alpha kx + k(\beta + k)y = e \cdot \dfrac{I}{d_i} - \delta k + (2\alpha\varepsilon - \beta\delta) \end{cases} \quad (11)$$

We remark here that the integer $I$ has exactly $2N$ integer divisors: $1 = d_1, d_2, \ldots, d_N = |I|$; and their negative counterparts $-d_1, -d_2, \ldots, -d_N$. When one makes a choice for one of the two factors on the left-hand side of (10); one chooses that factor to be $e \cdot d_i$, where $e = 1$ or $-1$. Then the other factor must equal $e \cdot \dfrac{I}{d_i}$; so that the product of the two factors equals $I$.

Below we use the well-known Kramer's rule from linear algebra in order to solve each linear system in (11). First, we compute the determinant of the matrix of the coefficients; which must be nonzero in order to apply Kramer's. We have,

$$\det\begin{bmatrix} 2\alpha k & k(\beta - k) \\ 2\alpha k & k(\beta + k) \end{bmatrix} = 2\alpha k \cdot (k(\beta + k)) - 2\alpha k(k(\beta - k))$$

$$= 4\alpha k^3 \neq 0, \text{ since } \alpha \neq 0 \text{ and k is a positive integer.}$$

Thus, each linear system in (11) will have a unique solution. Since all the coefficients (including the constant terms) are integers, it follows each of the $2N$ system will have a unique solution $(x, y)$ in which both $x$ and $y$ are rational numbers. Let $(x_i, y_i)$ be either of the two solutions (one for each system) of the two systems in (11); one system is obtained for $e = 1$, the other for $e = -1$. If we have to specify which one of the two solutions we are referring to; then $(x'_i, y'_i)$ will stand for the solution of the system in (11) with $e = 1$; while $(x''_i, y''_i)$ will be the solution of the system in (11) with $e = -1$. To be able to write down explicit formulas for $x_i$ and $y_i$; we must compute two more determinants as required by Kramer's rule. This is standard material that can be found not only in linear algebra texts, but also in college algebra and precalculus texts. Here is the end result after some simplifying:

$$\begin{cases} x_i = \dfrac{e \cdot d_i(\beta + k) - e \cdot \dfrac{I}{d_i}(\beta - k) - 2\delta k^2 - 2\beta(2\alpha\varepsilon - \beta\delta)}{4\alpha k^2} \\ y_i = \dfrac{e \cdot \dfrac{I}{d_i} - e \cdot d_i + 2(2\alpha\varepsilon - \beta\delta)}{2k^2} \end{cases} \quad (12)$$

## 4. A numerical example

Consider the hyperbola with equation
$$2x^2 - 5xy + 2y^2 - x + y - 1 = 0 \tag{13}$$

We have $\alpha = 2, \beta = -5, \gamma = 2, \delta = -1, \varepsilon = 1, J = -1$.
And $\beta^2 - 4\alpha\gamma = 25 - 4(2)(2) = 9 = k^2;\ k = 3$

Furthermore we obtain,
$2\alpha k = 12, k(\beta - k) = -24, \delta k = -3, 2\alpha\varepsilon - \beta\delta = -1,$
$k(\beta + k) = -6, \delta^2 - 4\alpha J = 9;$

and thus by (10a),
we also get $I = 80$. The integer 80 has exactly $N = 10$ positive divisors; and so the number of linear of systems to be solved is $2N = 20$. However, if we look at equation (10); in the case of this example we have,

$$(12x - 24y - 4)(12x - 6y - 2) = 80 \tag{14}$$

which shows that we can cancel out a common factor from each of the two sides of equation (13); and still have two factors on the left-hand side of the resulting equation, each of which is a linear polynomial in $x$ and $y$ with integer coefficients.
Specifically equation (14) is equivalent to

$$(3x - 6y - 1)(6x - 3y - 1) = 10 \tag{15}$$

So this observation about the common factor of 8 simply expedites the solving process. Equation (15) yields eight linear systems:

(a) $\begin{cases} 3x - 6y - 1 = 1 \\ 6x - 3y - 1 = 10 \end{cases}$  $\begin{cases} 3x - 6y - 1 = -1 \\ 6x - 3y - 1 = -10 \end{cases}$ (b)

(c) $\begin{cases} 3x - 6y - 1 = 2 \\ 6x - 3y - 1 = 5 \end{cases}$  $\begin{cases} 3x - 6y - 1 = -2 \\ 6x - 3y - 1 = -5 \end{cases}$ (d)

(e) $\begin{cases} 3x - 6y - 1 = 10 \\ 6x - 3y - 1 = 1 \end{cases}$  $\begin{cases} 3x - 6y - 1 = -10 \\ 6x - 3y - 1 = -1 \end{cases}$ (f)

(g) $\begin{cases} 3x - 6y - 1 = 5 \\ 6x - 3y - 1 = 2 \end{cases}$  $\begin{cases} 3x - 6y - 1 = -5 \\ 6x - 3y - 1 = -2 \end{cases}$ (h)

Of these 8 systems only (b),(c),(f), and (g) have integer solutions. These are also the solutions of equation (13):

$$(x, y) = (-2, -1), (1, 0), (1, 2), (0, -1).$$

If we had proceeded without cancelling the factoring the common factor 8 in equation (14); we would have found out that some of the resulting 20 systems are in fact equivalent; in the end, of course, we would have found the same four integer solutions above. Also note that the 10 positive integer divisors of 80 are:

$$d_1 = 1, d_2 = 2, d_3 = 4, d_4 = 5, d_5 = 8, d_6 = 10, d_7 = 16, d_8 = 20, d_9 = 40, \text{ and } d_{10} = 80.$$

5. **Observations and remarks**

   Keep in mind that among the $2N$ linear systems one must solve in order to find all the integer solutions of equation (10) and thus of (1)); there may a few or several groups, with each such group containing equivalent systems. When that happens, the whole solving process is significantly simplified/reduced.
   Now, let us take another look at the integer

$$I = k^2(\delta^2 - 4\alpha J) - (2\alpha\varepsilon - \beta\delta)^2.$$

If $\beta$ or $\delta$ is even; then by inspection we see that $I \equiv 0 \pmod{4}$. This statement is obvious when $\delta$ is even. When $\beta$ is even, then by $\beta^2 - 4\alpha\gamma = k^2$, it follows that $k$ is even as well; so must $I$.
If both $\beta$ and $\delta$ are odd integers; then so is $k$. And so, $k^2 \cdot \delta^2 \equiv 1 \equiv (2\alpha\varepsilon - \beta\delta)^2 \pmod{4}$, from which it easily follows that

$$I \equiv 1 - 1 \equiv 0 \pmod{4} \begin{pmatrix} \text{recall that the square of any odd integer, is in fact congruent to } 1 \pmod{8}; \\ \text{and thus in particular to } 1 \pmod{4} \end{pmatrix}$$

*Conclusion: The integer $I$ is always a multiple of 4*

6.    **The case $\beta = \delta = \varepsilon = 0$.**

When $\beta = \delta = \varepsilon = 0$, equation (1) becomes
$$\alpha x^2 + \gamma y^2 = -J \tag{16}$$

According to (2), we also have

$$-4\alpha\gamma = k^2 \tag{16a}$$

Equation (16a) says that the product $\alpha\gamma$ must equal minus an integer or perfect square. That can of course occur without one of $\alpha$ and $\gamma$ being an integer square while the other being minus a perfect square. For example $\alpha = 18$ and $\gamma = -50$, gives $\alpha\gamma = -2^2 \cdot 3^2 \cdot 5^2 = -(30)^2$. However, if one of $\alpha$ and $\gamma$ is an integer square while the other is minus a square; then this would be a sufficient condition while implies that $\alpha\gamma$ is minus an integer square. Assume then that,

$$\alpha = l^2 \text{ and } \gamma = -m^2 \qquad (16b)$$

Where $l$ and $m$ are positive integers.

And so, the value of the natural number $k$ is $k = 2ml$

It is clear from (16) that if $(x_0, y_0)$ is an integer solution, then so are the pairs $(x_0, -y_0), (-x_0, y_0)$ and $(-x_0, -y_0)$.

By (16) and (16a) we get,

$$\begin{cases} l^2 x^2 - m^2 y^2 = -J \\ or\ equivalently, (lx - my)(lx + my) = -J \end{cases} \qquad (17)$$

Note that if $J$ is an integer congruent to 2 modulo 4; then (17) has no integer solutions. This is true because $l^2 x^2 - m^2 y^2 = (lx)^2 - (my)^2 \equiv 0, 1, or\ 3 \pmod{4}$ according to whether both $lx$ and $my$ are even or odd; $lx$ is odd and $my$ even; or $lx$ is even and $my$ odd.

Below we examine the two cases when $J = -1, -p$ where $p$ is an odd prime; and so $-J = 1, p$

1. Suppose $-J = 1$. Since the only positive divisor of 1 is 1. To find all the integer solutions of (17), we must solve the two linear systems:

   (i) $\begin{cases} lx - my = 1 \\ lx + my = 1 \end{cases}$  $\qquad$ $\begin{cases} lx - my = -1 \\ lx + my = -1 \end{cases}$ (ii)

   The solution of $(i)$ is $(x, y) = \left(\frac{1}{l}, 0\right)$ and of $(ii)$ is $(x, y) = \left(-\frac{1}{l}, 0\right)$

   *We see that when $J = -1$, equation (17) will have exactly two integer solutions when $l = 1$, these being $(1, 0)$ and $(-1, 0)$. Otherwise, for any other value of the natural number $l$; and for positive integer $m$, equation (17) has no integer solutions.*

2. Suppose $-J = p$. The positive divisors of $p$ are $1$ and $p$. To find all the integer solutions of (17), we must solve the four linear systems below:

$$(iii) \begin{cases} lx - my = 1 \\ lx + my = p \end{cases} \quad (iv) \begin{cases} lx - my = p \\ lx + my = 1 \end{cases} \quad (v) \begin{cases} lx - my = -1 \\ lx + my = -p \end{cases}$$

$$(vi) \begin{cases} lx - my = -p \\ lx + my = -1 \end{cases}$$

The solutions of $(iii), (iv), (v)$, and $(vi)$ are respectively,

$$(x, y) = \left( \frac{p+1}{2l}, \frac{p-1}{2m} \right), \left( \frac{p+1}{2l}, \frac{1-p}{2m} \right), \left( \frac{-(p+1)}{2l}, \frac{1-p}{2m} \right) \text{ and } \left( \frac{-(p+1)}{2l}, \frac{p-1}{2m} \right).$$

We conclude that when $J = -p$, $p$ and odd prime; and $2l$ is not a divisor of $(p+1)$; (nonexclusive) or $2m$ is not a divisor of $(p-1)$; then equation (17) has no integer solutions. On the other hand, when $2l$ is a divisor of $(p+1)$ and $2m$ is a divisor $(p-1)$, then equation (17) has exactly four integer solutions (which are listed above)

7. **The case $\alpha = 1 = k$**

By (2) and $(10_a)$ we have

$$\begin{cases} \beta^2 - 4\gamma = 1 \\ \text{and } I = \delta^2 - 4J - (2\varepsilon - \beta\delta)^2 \end{cases} \tag{18}$$

Note that $\beta$ must be odd in this case, and so $\beta^2 \equiv 1 (\mod 8)$; and thus $\gamma$ even. Going back to (12),

$$\begin{cases} x_i = \dfrac{ed_i(\beta+1) - e\dfrac{I}{d_i}(\beta-1) - 2\delta - 2\beta(2\varepsilon - \beta\delta)}{4} \\ \\ y_i = \dfrac{e\dfrac{I}{d_i} - ed_i + 2(2\varepsilon - \beta\delta)}{2} \end{cases} \tag{19}$$

with $e = 1 \text{ or } -1$

As we know from the previous section, $I$ is always a multiple of 4. Consider the special case $I = 2^n$, $n \geq 2$.

The $N = n+1$ positive divisors of $I$ are,

$d_1 = 1, d_2 = 2, \ldots, d_n = 2^{n-1}, d_{n+1} = 2^n = d_N; d_i = 2^{i-1}$, for $i = 1 = 2, \ldots, n+1$

Note that when $i = 1; d_1 = 1, d_1 = 1, \dfrac{I}{d_1} = 2^n$ and therefore $y_1$ is half an odd integer; thus, $y_1$ is not an integer. Likewise when $i = n+1, d_{n+1} = 2^n$, but $\dfrac{I}{d_{n+1}} = 1$; again, $y_{n+1}$ is not an integer. On the other hand, for $2 \le i \le n-1$, both $x_i$ and $y_i$ are integers. Indeed in this case $d_i = 2^{i-1}$ and $\dfrac{I}{d_i} = 2^{n-i+1}$ are both even, and therefore $y_i$ is even by inspection. Also, since $\beta$ is odd, both $(\beta - 1)$ and $(\beta + 1)$ are even integers. Consequently, $ed_i(\beta + 1) \equiv e \cdot \dfrac{I}{d_i}(\beta - 1) \equiv 0 \pmod 4$. Furthermore, $2\delta \equiv 2\beta(2\varepsilon - \beta\delta) \pmod 4$, since $\beta \equiv 1 \pmod 2$; and so $\beta^2 \equiv 1 \pmod 4$. When $\delta$ is odd, then $2\delta \equiv 2\beta(2\varepsilon - \beta\delta) \equiv 2 \pmod 4$; while when $\delta$ is even $2\delta \equiv 2\beta(2\varepsilon - \beta\delta) \equiv 0 \pmod 4$. We see that when $I = 2^n, n \ge 2$, and the conditions in (18) are satisfied, the hyperbola described by equation (1) contains exactly $2(n-1)$ (distinct, see remark below) integral points. Using $d_i = 2^{i-1}$ and $\dfrac{I}{d_i} = 2^{n-i+1}$, the formulas in (19), and after some simplifying, we can state the following theorem.

**Theorem 1**
Suppose that $\alpha, \beta, \gamma, \delta, \varepsilon, J$ are integers such that
$\alpha = 1, \beta^2 - 4\gamma = 1, \gamma \ne 0,$ and $I = \delta^2 - 4J - (2\varepsilon - \beta\delta)^2 = 2^n, n \ge 2$. Then, the hyperbola described by equation (1) contains precisely $2(n-1)$ integral points given by the formulas

$$\begin{cases} x_i = \dfrac{e2^{i-2}(\beta + 1) - e \cdot 2^{n-i}(\beta - 1) - \delta - (2\varepsilon - \beta\delta)\beta}{2} \\ y_i = e2^{n-i} - e \cdot 2^{i-2} + 2\varepsilon - \beta\delta \end{cases}$$

for $i = 2, \ldots, n-1$; and with $e = 1$ or $-1$

**Remark 1.** In the last section, Section 10, we prove that the $2(n-1)$ integral points of Theorem 1, are in fact distinct.

8. **The case $I = 0$**

When $I = k^2(\delta^2 - 4\alpha J) - (2\alpha\varepsilon - \beta\delta)^2 = 0$, then the curve described by (1) and (2), is a pair of two intersecting straight lines. In this case, the solution sets of (10) and therefore of (1) as well is the union of two sets:

$S = S_1 \cup S_2$, where $S_1$ is the solution set of the linear diophantine equation,
$2\alpha k x + k(\beta - k)y + \delta k + 2\alpha\varepsilon - \beta\delta = 0$; or equivalently

$$2\alpha k x + k(\beta - k)y = -\delta k - (2\alpha\varepsilon - \beta\delta) \tag{20a}$$

And $S_2$ is the solution set of the linear diophantine equation

$$2\alpha k x + k(\beta + k)y = -\delta k + (2\alpha\varepsilon - \beta\delta) \tag{20b}$$

All four possibilities can occur: $S_1 = \phi$ (empty set) and $S_2 \neq \phi$; $S_1 \neq \phi$ and $S_2 = \phi$; $S_1 = S_2 = \phi$; or $S_1 \neq \phi$ and $S_2 \neq \phi$. In both cases of (20a) and (20b); we are dealing with a linear diophantine equation in two variables:

$ax + by = c$; with $a, b, c$ being integers.

This is well known, standard material, that can be found in almost every introductory number theory book. Such an equation has either no integer solutions or infinitely many solutions. It has infinitely many integer solutions if, and only if, the greatest common divisor $\gcd(a,b)$ of the coefficients $a$ and $b$; is also a divisor of $c$. There is a well known method for finding all the integer solutions of a two-variable linear diophantine equation with integer coefficients; a procedure based on what is known in the literature as the Euclidean algorithm for finding the greatest common divisor of two integers. In the end, all the integer solutions can be parametrically expressed in terms of one integer parameter. The reader may refer to [2] or [3] for further details. In the case of the linear diophantine equation (20a); equation (20a) will have infinitely many integer solutions if, and only if, $\gcd(2\alpha k, k(\beta - k)) = k \gcd(2\alpha, \beta - k)$ is a divisor of the integer $-\delta k - (2\alpha\varepsilon - \beta\delta)$. If that is the case, the set $S_1$ will contain infinitely many integral points; otherwise $S_1 \neq \phi$. Likewise, equation (20b) will have infinitely integer solutions exactly when $\gcd(2\alpha k, k(\beta + k)) = k \gcd(2\alpha, \beta + k)$ is a divisor of the integer $-\delta k + (2\alpha\varepsilon - \beta)$. If this is the case, the set $S_2$ will contain infinitely many integral points. Otherwise $S_2 \neq \phi$.

9. **The case $\delta = \varepsilon = J = 0$**

It immediately follows from (10a) that $I = 0$. Also note that in this case, equation (1) is a homogeneous quadratic equation in $x$ and $y$ (i.e. every monomial term has degree 2). Since $I = 0$, equations (20a) and (20b) from the previous section, do apply.

we obtain $\begin{cases} 2\alpha k x + k(\beta + k)y = 0 \\ \text{or equivalently, } 2\alpha x + (\beta + k)y = 0 \end{cases}$ (21a)

And also (20b) becomes
$$\begin{cases} 2\alpha k x + k(\beta - k)y = 0 \\ \text{or equivalently, } 2\alpha x + (\beta - k)y = 0 \end{cases} \quad (21b)$$

In this special case, the conic described by equations (1) and (2); is a pair of two intersecting straight lines, with their point of intersection being the origin $(0,0)$.

Below we find the sets $S_1$ and $S_2$. First $S_1$: we find all the integer solutions to (21a)
Let $d_1 = \gcd(2\alpha, \beta + k)$. Then $2\alpha = d_1 \cdot \rho_1$ and $\beta + k = d_1 \cdot v_1$, where $v_1$ and $\rho_1$ are relatively prime integers; $\gcd(v_1, \rho_1) = 1$.
Accordingly, by (21a) we obtain,

$$\rho_1 \cdot x = -v_1 y \quad (22a)$$

Since $\rho_1$ is relatively prime to $v_1$; and thus to $-v_1$ as well; and it divides the product $-v_1 y$; it must be a divisor of $y$ (this is the well known Euclid's Lemma; which is fundamental in proving the Fundamental Theorem (Unique Factorization) of Arithmetic) We set $y = \rho_1 \cdot t$; $t$ an integer.
Therefore from (22a)
we obtain $x = -v_1 t$. Thus,

$S_1 = \{(x, y) | x = -v_1 t, y = \rho_1 t, t \in Z\}$; where $\rho_1$ and $v_1$ are the integers defined by:

$$\rho_1 = \frac{2\alpha}{\gcd(2\alpha, \beta + k)}, \text{ and } v_1 = \frac{\beta + k}{\gcd(2\alpha, \beta + k)}$$

Next, if we let $d_2 = \gcd(2\alpha, \beta - k)$ and we set $2\alpha = d_2 \cdot \rho_2$ and $\beta - k = d_2 \cdot v_2$; $\rho_2, v_2 \in Z$ (21b)
with
$\gcd(\rho_2, v_2) = 1$. We obtain from the equation

$$\rho_2 \cdot x = -v_2 y \qquad (22b)$$

Working similarly as with the previous case, we find that

$$S_2 = \{(x, y) | x = -v_2 t, y = \rho_2 t, t \in Z\}$$

where $\rho_2 = \dfrac{2\alpha}{\gcd(2\alpha, \beta - k)}$, $v_2 = \dfrac{\beta - k}{\gcd(2\alpha, \beta - k)}$.

## 10. The distinctness of the $2(n-1)$ integral points in Theorem 1

In this section we prove that the $2(n-1)$ integral points, as defined in Theorem 1, are distinct. We have,

$$\begin{cases} x_i = \dfrac{e2^{i-2}(\beta+1) - e2^{n-i}(\beta-1) - \delta - \beta(2\varepsilon - \beta\delta)}{2} \\ y_i = e2^{n-i} - e2^{i-2} + 2\varepsilon - \beta\delta \end{cases} \qquad (23)$$

Recall that $\beta^2 - 4\gamma = 1$, and so $\beta$ is odd.

First we show that for each $i = 2, \ldots, n-1 \,(n \geq 2)$, the two pairs $(x_i, y_i)$ are distinct. One of these two pairs is obtained when $e = 1$; we will denote it by $(x'_i, y'_i)$; the other corresponds to $e = -1$, we will denote it by $(x''_i, y''_i)$. Let us compare their $y$-coordinates. When can it be $y'_i = y''_i$? We have
$$y'_i = y''_i \Leftrightarrow 2^{n-i} - 2^{i-2} + 2\varepsilon - \beta\delta = -2^{n-i} + 2^{i-2} + 2\varepsilon - \beta\delta \Leftrightarrow 2^{n-i+1} = 2^{i-1} \Leftrightarrow n - i + 1 = i - 1$$
$$\Leftrightarrow i = \dfrac{n+2}{2},$$

which is possible only when $n$ is even. But then, for $i = \dfrac{n+2}{2}$, can the $x$-coordinates be equal?

We have
$$\begin{cases} x'_i = x''_i \\ i = \dfrac{n+2}{2} \end{cases} \Leftrightarrow \begin{cases} 2^{i-2}(\beta+1) - 2^{n-i}(\beta-1) = -2^{i-2}(\beta+1) + 2^{n-i}(\beta-1) \\ i = \dfrac{n+2}{2} \end{cases}$$

$$\Leftrightarrow \begin{cases} 2^{i-1}(\beta+1) = 2^{n-i+1}(\beta-1) \\ i = \dfrac{n+2}{2} \end{cases} \Leftrightarrow 2^{\frac{n}{2}}(\beta+1) = 2^{\frac{n}{2}}(\beta-1).$$

$\Leftrightarrow \beta + 1 = \beta - 1 \Leftrightarrow 2 = 0$, an obvious impossibility.

Next, we prove that no two integral points $(x_i, y_i)$ and $(x_j, y_j)$, can coincide when $i \neq j$. Assume then that $2 \leq j < i \leq n-1$. In particular $i - j \geq 1$.

For the point $(x_i, y_i)$ we will use $e_i$ in the formulas (23); where $e_i = 1\, or -1$. Likewise for the point $(x_j, y_j)$, we will use $e_j$ in the formulas (23); where $e_j = 1\, or -1$.

First, let us examine the possibility $y_i = y_j$. This is equivalent to

$$e_i 2^{n-i} - e_i 2^{i-2} = e_j 2^{n-j} - e_j 2^{j-2} \Leftrightarrow 2^{n-i} \cdot [e_i - e_j 2^{i-j}] = 2^{j-2} \cdot [e_i 2^{i-j} - e_j] \tag{24}$$

Since $i - j \geq 1$; both integers $(e_i - e_j 2^{i-j})$ and $(e_i \cdot 2^{i-j} - e_j)$ are odd; therefore equation(24) implies the two powers of 2 on either side of (24) must be equal; which means that we must have

$$\begin{aligned} n - i &= j - 2; \\ i + j &= n + 2 \end{aligned} \tag{25}$$

Accordingly (24) yields,

$$\begin{aligned} e_i - e_j 2^{i-j} &= e_i 2^{i-j} - e_j \\ \Leftrightarrow (e_i + e_j)(2^{i-j} - 1) &= 0 \end{aligned} \tag{26}$$

However $i - j \geq 1$; and so $2^{i-j} - 1 > 0$.
Thus (26) implies $e_i + e_j = 0$; which means that

$$\begin{cases} either\ e_i = 1\, and\, e_j = -1 \\ or\ e_i = -1\, and\, e_j = 1 \end{cases} \tag{27}$$

We have seen that $y_i = y_j$, precisely when the conditions (25) and (27) are satisfied. Now, can in addition to this, also have $x_i = x_j$?

Going back to (23), we see that $x_i = x_j$ is equivalent to

$$e_i 2^{i-2}(\beta + 1) - e_j 2^{n-i}(\beta - 1) = e_j 2^{j-2}(\beta + 1) - e_j 2^{n-j}(\beta - 1)$$
$$\Leftrightarrow (\beta + 1)[e_i 2^{i-2} - e_j 2^{j-2}] = (\beta - 1)[e_i 2^{n-i} - e_j 2^{n-j}]$$

From (25) we know that $j = n + 2 - i$. Substituting for $j$ in the last equation yields

$$(\beta+1)\left[e_i 2^{i-2} - e_j 2^{n-i}\right] = (\beta-1)\left[e_i 2^{n-i} - e_j 2^{i-2}\right]$$

Also, by (27), we have $e_j = -e_i$ Thus, the last equation implies

$$(\beta+1)e_i\left[2^{i-2} + 2^{n-i}\right] = (\beta-1)e_i\left[2^{n-i} + 2^{i-2}\right];$$

which gives $\beta + 1 = \beta - 1$, an impossibility.

This proves that $(x_i y_i) \neq (x_j y_j)$, for all $i$ and $j$ such that $2 \leq j < i \leq n - 1$. We are done    □